\documentclass[journal,10pt,onecolumn,draftclsnofoot]{IEEEtran}

\usepackage{graphicx}
\usepackage{float}
\ifCLASSINFOpdf
\else
\fi

\usepackage{amssymb, amscd, amsthm}
\usepackage{amsmath}
\allowdisplaybreaks
\usepackage{color}
\usepackage[ruled,norelsize]{algorithm2e}
\makeatletter
\newcommand{\removelatexerror}{\let\@latex@error\@gobble}
\makeatother

\ifCLASSOPTIONcompsoc
    \usepackage[caption=false, font=normalsize, labelfont=sf, textfont=sf]{subfig}
\else
\usepackage[caption=false, font=footnotesize]{subfig}
\fi

\newtheorem{thm}{Theorem}[section]

\newtheorem{prop}{Proposition}[section]

\usepackage{url}
\usepackage{cite}
\usepackage{bm}
\usepackage{multirow}
\renewcommand{\baselinestretch}{.97} 
\begin{document}
%
\title{The Value of Including Unimodality Information in Distributionally
Robust Optimal Power Flow}

\author{Bowen~Li,~\IEEEmembership{Member,~IEEE,}
        Ruiwei~Jiang,~\IEEEmembership{Member,~IEEE,}
        and~Johanna~L. Mathieu,~\IEEEmembership{Senior Member,~IEEE}%
\thanks{This research was supported by the U.S. National Science Foundation Awards CCF-1442495 and CMMI-1662774. B.~Li is with Argonne National Laboratory, Lemont, IL, USA. email: bowen.li@anl.gov.  J.~Mathieu is with the Department of Electrical Engineering \& Computer Science, University of Michigan, Ann Arbor, MI, USA. email: jlmath@umich.edu. R.~Jiang is with the Department of Industrial $\&$ Operations Engineering, University of Michigan, Ann Arbor, MI, USA. email: ruiwei@umich.edu.}}%


\maketitle

\begin{abstract}
 To manage renewable generation and load consumption uncertainty, chance-constrained optimal power flow (OPF) formulations and various solution methodologies have been proposed. However, conventional solution approaches often rely on accurate estimates of uncertainty distributions, which may not exist. When the distributions are not known but can be limited to a set of plausible distributions, termed an ambiguity set,
 distributionally robust (DR) optimization can be used to ensure that chance constraints hold for all distributions in that set. However, DR OPF yields conservative solutions if the ambiguity set is too large.  In this paper, we assess the value of using both moment and unimodality information, which shrinks the ambiguity set and reduces conservatism, in DR OPF problems. Most practical uncertainty distributions in power systems are unimodal. Exact reformulations, approximations, and efficient solving techniques were developed in a previous paper. This paper develops an optimal parameter selection approach that searches for an optimal approximation, significantly improving the computational efficiency and solution quality. We evaluate the performance of the approach against existing chance-constrained OPF approaches using modified IEEE 118-bus and 300-bus systems with high penetrations of renewable generation. Results show that including unimodality information reduces solution conservatism and cost without significantly degrading reliability.

\end{abstract}

\begin{IEEEkeywords}
Optimal power flow, chance constraint, distributionally robust optimization, $\alpha$-unimodality
\end{IEEEkeywords}

\IEEEpeerreviewmaketitle

\section{Introduction} \label{sec: intro}

\IEEEPARstart{P}{revious} research has developed approaches to ensure power system reliability under uncertainties (such as renewable generation forecast error) by chance-constrained optimal power flow (OPF) models, in which physical constraints are required to be satisfied with high probability, e.g., \cite{cc1,cc2,vrakopoulou_probabilistic_2013,bienstock_chance_2014,cc4,mariatsg,bowentsg}. Conventional approaches to solving chance-constrained OPF problems include scenario approximation \cite{Campi2009, Margellos2014}, analytical reformulations based on known distributions (e.g., Gaussian) \cite{bienstock_chance_2014,cc4,bowentsg,roald1}, and sample average approximation (SAA) \cite{SAA1,SAA2}. Scenario approximation approaches rely on a large number of scenarios and often provide overly-conservative results.
Analytical reformulations usually require less computational effort; however, it is often difficult to accurately estimate the \emph{joint} probability distribution of the uncertain parameters and so solutions can be unreliable. SAA performs better as the number of samples increases, but that also increases its computational burden as more binary variables and constraints are needed when recasting the SAA formulation as a mixed-integer program. 

In contrast, distributionally robust (DR) optimization ensures that chance constraints hold with regard to all probability distributions within an ambiguity set \cite{el2003worst,delage2010distributionally,stellato2014data,ruiwei}. The approach is closely related to both robust and stochastic optimization because 1) it reduces to robust optimization if the ambiguity set includes only the support information and 2) it reduces to a  chance-constrained program if the ambiguity set includes only a single distribution. By incorporating distributional information of the uncertainty (such as moments) into the ambiguity set, DR optimization can achieve a better trade-off between solution costs and reliability than the aforementioned existing approaches. The conservatism of the DR approach is related to the ambiguity set: if it includes unrealistic distributions, then the solution may be more costly than necessary. A recent thrust of research in DR optimization is the development of methods to incorporate additional information, e.g., on the distribution structure, into the ambiguity set so that unrealistic distributions can be eliminated. However, incorporating additional information often comes with additional computational burden.

The objective of this work is to assess the value of using both moment and structural information, specifically, unimodality, in the DR OPF model. We investigate the trade-off between solution quality (cost, reliability) and computational efficiency, and compare our approach with a variety of existing ones. Our goal is to determine whether the additional computational burden is worth the improvement in solution quality. 

Previous DR OPF research has derived tractable reformulations assuming ambiguity sets based on moments \cite{dr5,summers,dr6,DRyiling,dr3,dr9_moment_2018,dr10_moment_2018}, discrepancy measure \cite{dr4,dr7,dr8}, and structural information such as symmetry \cite{roald1}, unimodality \cite{roald1,DRbowen,UnimodalBL,summers}, and log-concavity \cite{bowenpscc}. Reference \cite{dr3} considered two-sided joint chance constraints for generator and transmission line limits and \cite{dr4,dr8,dr7} constructed the ambiguity set based on the discrepancy between the real distribution and the empirical distribution. Here, we consider an ambiguity set that incorporates the first two moments and a generalized $\alpha$-unimodality \cite{unimodal}, which is typically satisfied by uncertainties in OPF models, such as wind power forecast error. Our prior work \cite{UnimodalBL,DRbowen} developed exact reformulations, approximations, and efficient solving techniques that we leverage here.

In this paper, we develop an optimal parameter selection (OPS) approach that helps us construct a high-quality conservative approximation of the DR chance constraints, which significantly improves upon the approximations in \cite{UnimodalBL,DRbowen}. The main step in the OPS approach is to find the closest piece-wise linear (PWL) outer approximation of a concave function that is independent of the values of the decision variables. We investigate multiple online and offline options to construct the approximation. We also provide the mathematical proofs of optimality and existence of the PWL approximation, and a heuristic solving algorithm. Then, we compare the approach to existing ones through case studies on modified 118-bus and 300-bus systems with high wind power penetration.

The remainder of the paper is organized as follows. In Section~\ref{sec: main}, we review some fundamental concepts and generalize the DR formulations in \cite{UnimodalBL}. In Section~\ref{sec: param_select}, we derive the OPS approach. In Section~\ref{sec: case}, we compare the performance of the new approach to existing ones and discuss the value of including unimodality information in the DR OPF problem. Section~\ref{sec: conclusion} concludes the paper.

\section{Distributionally Robust Chance Constraints}  \label{sec: main}

In this section, we review DR chance constraints and generalize the results from \cite{UnimodalBL}.
We assume constraints with uncertainty can be transformed into the form
\begin{equation}
a(x)^{\top}\xi\leq b(x), \label{eq: cons}
\end{equation}
where $x\in\mathbb{R}^n$ represents decision variables and $a(x):\mathbb{R}^n\rightarrow\mathbb{R}^l$ and $b(x):\mathbb{R}^n\rightarrow\mathbb{R}$ represent two affine functions of $x$. Uncertainty $\xi\in\mathbb{R}^l$ is defined on probability space $(\mathbb{R}^l,\mathcal{B}^l,\mathbb{P}_{\xi})$ with Borel $\sigma$-algebra $\mathcal{B}^l$ and probability distribution $\mathbb{P}_{\xi}$. 
To ensure \eqref{eq: cons} is satisfied with at least a probability threshold $1-\epsilon$, we define the chance constraint
\begin{equation}
\mathbb{P}_{\xi}\left(a(x)^{\top}\xi\leq b(x)\right)\geq1-\epsilon,\label{eq: cc}
\end{equation}
where $1-\epsilon$ normally takes a large value (e.g., $0.99$ \cite{charnes1958cost,miller1965chance}).

We consider two types of ambiguity sets. 
The first includes moment information only
\begin{align}
\mathcal{D}_{\xi} := \left\{\mathbb{P}_{\xi} \in \mathcal{P}^l: \mathbb{E}_{\mathbb{P}_{\xi}}[\xi] = \mu, \ \mathbb{E}_{\mathbb{P}_{\xi}}[\xi \xi^{\top}] = \Sigma\right\}, \label{Dset1}
\end{align}
and the second includes moment and unimodality information
\begin{align}
&\mathcal{U}_{\xi} := \bigl\{\mathbb{P}_{\xi} \in \mathcal{P}^l_{\alpha} \cap \mathcal{D}_{\xi}: \ \mathcal{M}(\xi)=m \bigr\}, \label{Dset2}
\end{align}
where $\mathcal{P}^l_{\alpha}$ and $\mathcal{P}^l$ denote all probability distributions on $\mathbb{R}^l$ with and without the requirement of $\alpha$-unimodality respectively; $\mu$ and $\Sigma$ denote the first and second moments of $\xi$; and $\mathcal{M}(\xi)=m$ specifies that the true mode value of $\xi$ is $m$. The value of $\alpha$ determines the shape of the unimodal distribution \cite{unimodal}. When $\alpha=1$, all the marginal distributions are univariate unimodal (i.e., the density function has a single peak called the mode and decaying tails). When $\alpha=l$, the density function of $\xi$ has a single peak at the mode and is non-increasing along any rays from the mode. As $\alpha\rightarrow\infty$, the requirement of unimodality gradually relaxes until it disappears. In practice, most uncertainties such as wind power forecast error follow a ``bell-shaped" unimodal distribution. 





The DR chance constraint with ambiguity set $\mathcal{D}_{\xi}$ is
\begin{equation}
\inf_{\mathbb{P}_{\xi}\in \mathcal{D}_\xi}\mathbb{P}_{\xi}\left(a(x)^{\top}\xi\leq b(x)\right)\geq1-\epsilon. \label{eq: drcc_old}
\end{equation}
\begin{thm}{(Theorem 2.2 in \cite{wagner2008stochastic})}\label{thm: drcc_set1}
The DR chance constraint \eqref{eq: drcc_old} can be exactly reformulated as
\begin{equation}
\sqrt{\left(\frac{1-\epsilon}{\epsilon}\right)a(x)^{\top}(\Sigma-\mu\mu^{\top})a(x)}\leq b(x)-a(x)^{\top}\mu.
\end{equation}

\end{thm}

Reference \cite{UnimodalBL} derives an exact reformulation for the DR chance constraint with $\mathcal{U}_{\xi}$, i.e.,
\begin{equation}
\inf_{\mathbb{P}_{\xi}\in \mathcal{U}_\xi}\mathbb{P}_{\xi}\left(a(x)^{\top}\xi\leq b(x)\right)\geq1-\epsilon; \label{eq: drcc}
\end{equation}
however, the results are derived assuming the mode is at the origin. Without loss of generality, we can rewrite \eqref{eq: cons} as $a(x)^{\top}(\xi-m)\leq b(x)-a(x)^{\top}m$ with $\xi-m$ as our new random vector whose mode is at the origin and generalize the results from \cite{UnimodalBL} to the case in which the mode is not necessarily at the origin.

\begin{thm}{(Extension of Theorem 1 in \cite{UnimodalBL})}\label{thm: drcc_set2}
The DR chance constraint \eqref{eq: drcc} can be exactly reformulated as
\begin{align}
&\sqrt{\frac{1-\epsilon-\tau^{-\alpha}}{\epsilon}}\|\Lambda a(x)\|\leq \tau\left(b(x)-a(x)^{\top}m\right)\nonumber\\
&-\left(\frac{\alpha+1}{\alpha}\right)(\mu-m)^{\top}a(x),\quad\forall\tau\geq\left(\frac{1}{1-\epsilon}\right)^{1/\alpha},\label{eq: drcc_set2}
\end{align}
where $\Lambda:=\left(\left(\frac{\alpha+2}{\alpha}\right)(\Sigma-\mu\mu^{\top})-\frac{1}{\alpha^2}(\mu-m)(\mu-m)^{\top}\right)^{1/2}$.
\end{thm}
\vspace{.1cm}
Since parameter $\tau$ can take an infinite number of values, the reformulation in Theorem~\ref{thm: drcc_set2} also involves an infinite number of second-order conic (SOC) constraints. To solve an optimization problem with \eqref{eq: drcc_set2}, \cite{UnimodalBL} proposes the following algorithm. 
\vspace{-.3cm}
\begin{figure}[H]
 \removelatexerror
  \begin{algorithm}[H]
   \caption{Iterative solving algorithm \cite{UnimodalBL}}
   Initialization: $i=1$, $\tau_0=\bigl(\frac{1}{1-\epsilon}\bigr)^{1/\alpha}$\;
   \vspace{0.1cm}
   Iteration $i$:

   Step 1: Solve the reformulated optimization problem with \eqref{eq: drcc_set2} using $\tau_j$ for all $j=0,\ldots,i-1$ and obtain optimal solution $x_i^*$. All $\tau_j$ values are collected from previous iterations;

   Step 2 ({\bf Separation}): Find worst case $\tau^*$ that results in the largest violation of \eqref{eq: drcc_set2} under $x_i^*$: {\bf IF} $\tau^*$ does not exist, {\bf STOP} and {\bf RETURN} $x_i^*$ as optimal solution; {\bf ELSE} {\bf GOTO} Step 3\;

   Step 3: Set $\tau_i=\tau^*$ and $i=i+1$, {\bf GOTO} Step 1.
  \end{algorithm}
  \label{Fig: algorithm1}
\end{figure}
\vspace{-.3cm}
\noindent The optimization problem in Step 1 can be solved directly. To efficiently perform Step 2 in Algorithm 1, we follow Proposition 3 in \cite{UnimodalBL} modified to consider modes at $m$.

Reference \cite{UnimodalBL} also developed a
 sandwich approximation to bound the optimal objective cost from both below and above. The approximation is asymptotic in the sense that it converges to the optimal objective cost with more parameters included. 
\begin{prop}{Relaxed Approximation (Extension of Proposition 4 in \cite{UnimodalBL})}\label{prop: drcc_set2_lb}
For integer $K\geq 1$ and real numbers $\tau_0\leq n_1<n_2\ldots<n_K\leq\infty$, \eqref{eq: drcc} implies the SOC constraints
\begin{align}
&\sqrt{\frac{1-\epsilon-{n_k}^{-\alpha}}{\epsilon}}\|\Lambda a(x)\|\leq {n_k}\left(b(x)-a(x)^{\top}m\right)\nonumber\\
&-\left(\frac{\alpha+1}{\alpha}\right)(\mu-m)^{\top}a(x),\quad\forall k=1,\ldots,K.\label{eq: drcc_set2_lb}
\end{align}
\end{prop}
\begin{prop}{Conservative Approximation (Extension of Proposition 5 in \cite{UnimodalBL})}\label{prop: drcc_set2_ub}
For integer $K\geq 2$ and real numbers $\tau_0 = n_1<n_2\ldots<n_K=\infty$, define a PWL function containing $(K-1)$ pieces:
\begin{align}
&g(\tau)=\min_{k=2,\ldots,K}\Biggl\{\sqrt{\frac{1}{\epsilon(1-\epsilon-n_k^{-\alpha})}}\Biggl[\left(\frac{\alpha n_k^{-\alpha-1}}{2}\right)\tau\nonumber\\
&\hspace{3cm} +1-\epsilon-\left(1+\frac{\alpha}{2}\right)n_k^{-\alpha}\Biggr]\Biggr\}.
\end{align}
Set $q_1=\tau_0$ and denote $q_2<\ldots<q_{K-1}$ as the $(K-2)$ break points of function $g(\tau)$. Then, \eqref{eq: drcc} is implied by the SOC constraints
\begin{align}
&g(q_k)\|\Lambda a(x)\|\leq {q_k}\left(b(x)-a(x)^{\top}m\right)\nonumber\\
&-\left(\frac{\alpha+1}{\alpha}\right)(\mu-m)^{\top}a(x),\quad\forall k=1,\ldots,K-1.\label{eq: drcc_set2_ub}
\end{align}
\end{prop}

The convergence of the sandwich approximation as parameter dimension increases is directly affected by the choice of parameters $n_k$. In \cite{UnimodalBL}, we proposed an online parameter selection approach that sets $n_k=\tau^*$, where $\tau^*$ is the worst case parameter determined in Step 2 of the $k$-th iteration of Algorithm 1. However, this selection of $n_k$ is only critical to the relaxed approximation in the sense that it only requires the reformulation to be satisfied for a finite number of critical $\tau$ values within the infinite number of values (see Theorem \ref{thm: drcc_set2} and Proposition \ref{prop: drcc_set2_lb}). It does not have direct connections to the conservative approximation. In the next section, we propose a new parameter selection approach to improve the conservative approximation.



\section{Optimal Parameter Selection}  \label{sec: param_select}
In this section, we propose an OPS approach for the conservative approximation of the DR chance constraint with $\mathcal{U}_{\xi}$ \eqref{eq: drcc}. Based on \cite{UnimodalBL}, $q_k$ for $k = 2, ..., K-1$ in Proposition~\ref{prop: drcc_set2_ub} define the break points of a concave PWL function $g(\tau)$ that outer approximates the nonlinear function
\begin{align}
v(\tau)=\sqrt{\frac{1-\epsilon-{\tau}^{-\alpha}}{\epsilon}},\label{eq: ops}
\end{align}
where $\tau\in[\tau_0,\infty)$. The OPS problem finds the optimal PWL outer approximation of $v(\tau)$. Conventional approaches to finding optimal PWL approximations  \cite{cox1971pwl,Imamoto2008pwl,vandewalle1975pwl} are not applicable to our problem as they do not consider outer approximations and assume the function has a bounded domain. Therefore, we build upon this prior work to make the following contributions.
\begin{enumerate}
\item We provide and prove optimality conditions for the optimal PWL outer approximation of $v(\tau)$, and we prove existence of a function that satisfies these conditions.
\item We develop a heuristic algorithm to search for the optimal PWL outer approximation.
\end{enumerate}

\subsection{Optimality and Existence}\label{sec: param_proof}
First, we define an optimal PWL approximation for our problem. Denote an $|\mathcal{S}|$-piece PWL outer approximation of $v(\tau)$ as $h(\tau)=\min_{s\in\mathcal{S}}\{d_s\tau+f_s\}$, where $\mathcal{S}$ is the set of indices representing the pieces and $d_s$ is non-increasing with increasing $s$. Also, denote $h_s(\tau)$ as the $s$-th piece in $h(\tau)$ and its domain as $\mathcal{H}_s$. Then, the error $e^{\text{max}}$ of the PWL outer approximation can be defined as the largest distance between $v(\tau)$ and its approximation
\begin{equation}
e^{\text{max}}=\max_{s\in\mathcal{S}}\max_{\tau\in\mathcal{H}_s}\{d_s\tau+f_s-v(\tau)\}.\label{eq: largesterr}
\end{equation}
Define the optimal PWL outer approximation as the one that minimizes $e^{\text{max}}$.


\begin{thm}{(Optimality)}\label{thm: pwl_multi}
For $|\mathcal{S}|\geq1$, an $|\mathcal{S}|$-piece PWL function  $h(\tau)=\min_{s\in\mathcal{S}}\{d_s\tau+f_s\}$ is an optimal PWL outer approximation of $v(\tau)$ if the following three conditions hold.
\begin{enumerate}

\item $h_{|S|}(\tau) = \sqrt{\frac{1-\epsilon}{\epsilon}}$.

\item $h_s(\tau)$ is tangent to $v(\tau)$ for all $s\in\mathcal{S}$.

\item $h(B_{s_1})-v(B_{s_1})=h(B_{s_2})-v(B_{s_2}),\ \forall s_1,s_2\in\mathcal{S}$, where $\{B_s,\ s \in \mathcal{S}\}$ are the break points and left end point of $h(\tau)$.

\end{enumerate}

\end{thm}
\noindent The proof is given in Appendix~\ref{appd: proof1}. 

\begin{thm}{(Existence)}\label{thm: pwl_multi_exist}
There always exists an $|\mathcal{S}|$-piece PWL function $h(\tau)=\min_{s\in\mathcal{S}}\{d_s\tau+f_s\}$ that satisfies all three conditions in Theorem~\ref{thm: pwl_multi}. 

\end{thm}

\noindent The proof is given in Appendix~\ref{appd: proof2}. 


\subsection{Algorithm}\label{VII-sec: param_alg}
Here we provide a heuristic algorithm to search for the optimal $|\mathcal{S}|$-piece PWL approximation of $v(\tau)$ on $[\tau_0,\infty)$. The algorithm is adapted from the recursive descent algorithm in \cite{Imamoto2008pwl}. We first define the following notation.
\begin{itemize}

\item $\Delta^i\in\mathbb{R}^{\mathcal{I}}$: the step size in iteration $i$

\item $\delta$: the percentage tolerance for termination criteria

\item $\mathcal{I}$: the maximum number of iterations

\item $\hbar^i(\tau)$: the first $|\mathcal{S}|-1$ pieces of the approximation in iteration $i$; we exclude the last zero-slope piece since it is trivial; $\hbar^i_s(\tau)$ is the $s$-th piece of $\hbar^i(\tau)$

\item $B^i\in\mathbb{R}^{|\mathcal{S}|}$: the break points and end points of $\hbar^i(\tau)$ in iteration $i$; $B^i_s$ is the $s$-th entry, where $B_1^i=\tau_0$

\item $T^i\in\mathbb{R}^{|\mathcal{S}|-1}$: the points at which $\hbar(\tau)$ is tangent to $v(\tau)$  in iteration $i$; $T^i_s$ is the $s$-th entry

\item $E^i\in\mathbb{R}^{|\mathcal{S}|}$: the distances between $\hbar(\tau)$ and $v(\tau)$ at all $B^i_s$  in iteration $i$; $E^i_s$ is the $s$-th entry; define the distance between last zero-slope piece and $v(\tau)$ at $B^i_{|\mathcal{S}|}$ as $e^i_T=\sqrt{(1-\epsilon)/{\epsilon}}-v(B^i_{|\mathcal{S}|})$


\end{itemize}

\begin{figure}[ht]
 \removelatexerror
  \begin{algorithm}[H]
   \caption{Heuristic searching algorithm}
   Initialization: $i=1$, $\Delta^1=1$, $\delta=0.01$, $\mathcal{I}=50$, $B^1_{|\mathcal{S}|}=10$,  $T^1_s=\tau_0+\frac{s(B^1_{|\mathcal{S}|}-\tau_0)}{|\mathcal{S}|} $ for $s=1,...,|\mathcal{S}|-1$;
   
   \vspace{0.1cm}
   Iteration $i$:

   Step 1: {\bf IF} $i\leq \mathcal{I}$, calculate $B_s^i$ for $s=2,...,|\mathcal{S}|-1$ by solving  $\hbar^i_{s-1}(B_s^i)=\hbar^i_{s}(B_s^i)$, where $\hbar^i_{s}(\tau)=v^{\prime}(T^i_s)(\tau-T^i_s)+v(T^i_s)$; {\bf ELSE} {\bf STOP} and {\bf RETURN} no convergence under current initialization.

   Step 2: Calculate $E^i$ and $e^i_T$;
 {\bf IF} $(1+\delta)e^i_T<E^i_{|\mathcal{S}|}$ or $(1+\delta)E^i_{|\mathcal{S}|}<e^i_T$, set $B^i_{|\mathcal{S}|} = 0.5(\hat{\tau}+B^i_{|\mathcal{S}|})$, where $\hat{\tau}$ is the solution of $\hbar^i_{|\mathcal{S}|-1}(\hat{\tau})=\sqrt{(1-\epsilon)/{\epsilon}}$, and {\bf GOTO} Step 1; {\bf ELSE} {\bf GOTO} Step 3;

  Step 3: {\bf IF} $\max(E^i)\leq (1+\delta)\min(E^i)$, {\bf STOP} and {\bf RETURN} $\hbar^i(\tau)$ with last zero-slope piece as optimal solution; {\bf ELSE} {\bf GOTO} Step 4\;

  Step 4: {\bf IF} $i=1$, {\bf GOTO} Step 5; {\bf ELSEIF} $\max(E^i)>\max(E^{i-1})$, set $i=i-1$, $\Delta^i=\Delta^i/2$, and {\bf GOTO} Step 5; {\bf ELSE} {\bf GOTO} Step 5\;

  Step 5: For  $s=1,...,|\mathcal{S}|-1$,
 calculate $$T^{i+1}_s=T^i_s+\frac{\Delta^i(E^i_{s+1}-E^i_{s})}{\frac{E^i_{s+1}}{B^i_{s+1}-T^i_s}+\frac{E^i_{s}}{T^i_s-B^i_s}}$$
 and set 
 $i=i+1$ and {\bf GOTO} Step 1\;

  \end{algorithm}
  \label{Fig: algorithm2}
\end{figure}

\begin{figure}
\centering
\includegraphics[width=5in]{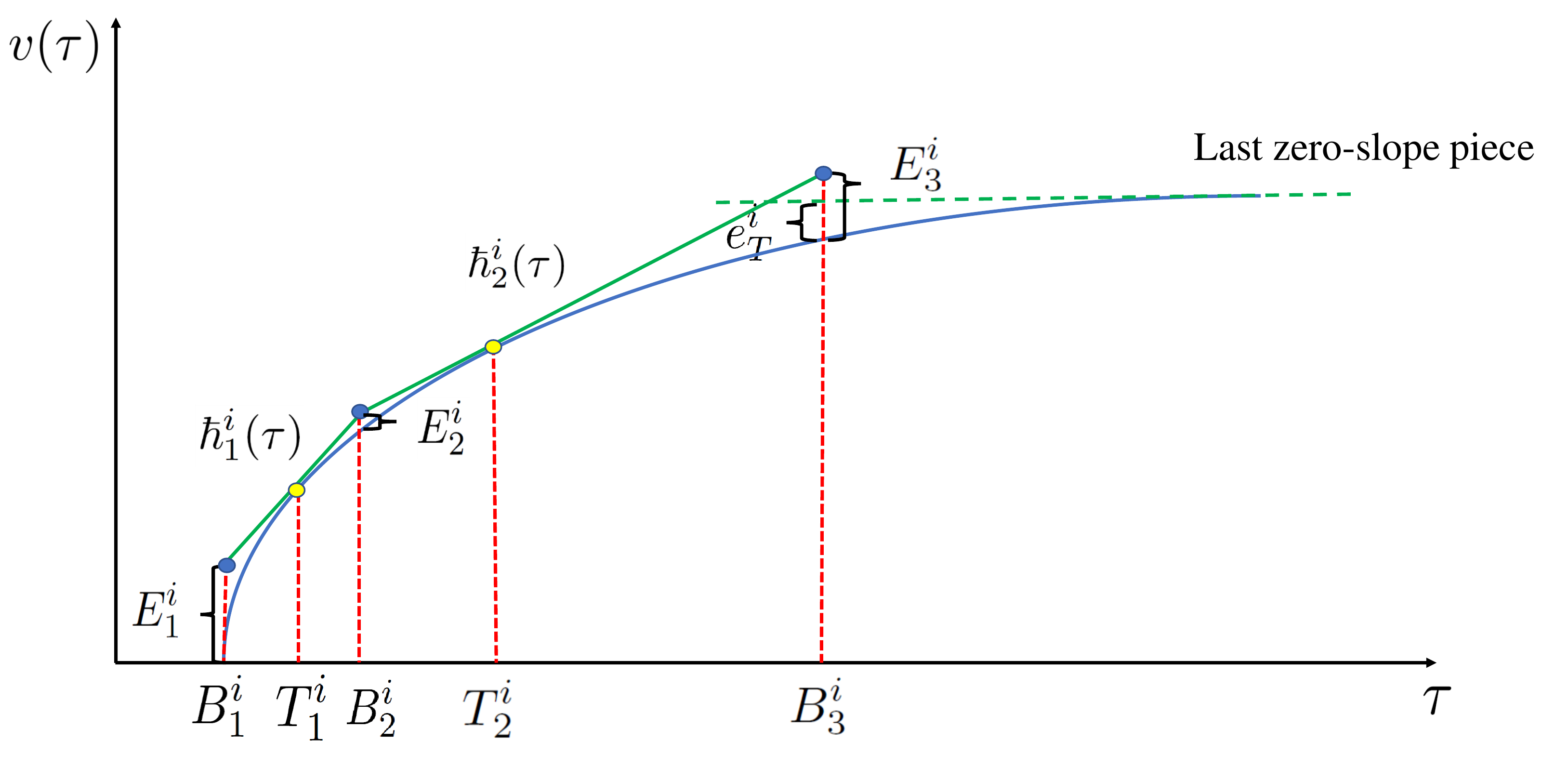}
\caption{Illustrative example of  Algorithm 2, iteration $i$ for $|\mathcal{S}|=3$. Function $h_i(\tau)$ (green solid lines) is a PWL outer approximation of $v(t)$ (blue curve) with tangent points $T^i$ (yellow dots) and break points and end points $B^i$ (blue dots). At convergence,  $E^{*}_1=E^{*}_2=E^{*}_3=e^{*}_T$.}
\label{fig: alg3_slides}
\end{figure}

The algorithm is given below (Algorithm 2). Figure~\ref{fig: alg3_slides} shows an illustrative example for $|\mathcal{S}|=3$. First, $\Delta^1$, $\delta$,  $\mathcal{I}$, and $B^1_{|\mathcal{S}|}$ are initialized (for different $v(\tau)$, they could take other reasonable values) and $T^1$ is computed by evenly dividing $[\tau_0,B^1_{|\mathcal{S}|}]$ into $|\mathcal{S}|$ segments.  In each iteration $i$, Step 1 calculates $B^i$ using $T^i$ and $B^i_{|\mathcal{S}|}$. Step 2 calculates $E^i$ and $e^i_T$, and coarsely adjusts $B^i_{|\mathcal{S}|}$. Since, 
at convergence ($i=*$) $E_s^* \approx e_T^* \ , \forall s\in\mathcal{S}$ with tolerance $\delta$, we should reduce $B^i_{|\mathcal{S}|}$ for $E^i_{|\mathcal{S}|}>e^i_T$, and vice versa. Step 3 checks for convergence.
Step 4 repeats the iteration with a smaller step size if the previous step did not produce an improvement. 
Step 5 adjusts $T^i$ to further reduce the differences among $E^i_s$. The adjustment is based on the approach in \cite{Imamoto2008pwl},
The optimal break points and end points $B^{*}$ can be used in Proposition~\ref{prop: drcc_set2_ub} to establish the corresponding conservative approximation.

\vspace{0.2cm}
\noindent{\bf Remark:} The optimal parameters obtained from Algorithm 2 are unique given a choice of $|\mathcal{S}|$. They are independent of the decision variables but dependent on the system parameters. Hence, they can be determined offline.

\subsection{Performance}

\begin{figure}
\centering
\vspace{.cm}
\includegraphics[width=5in]{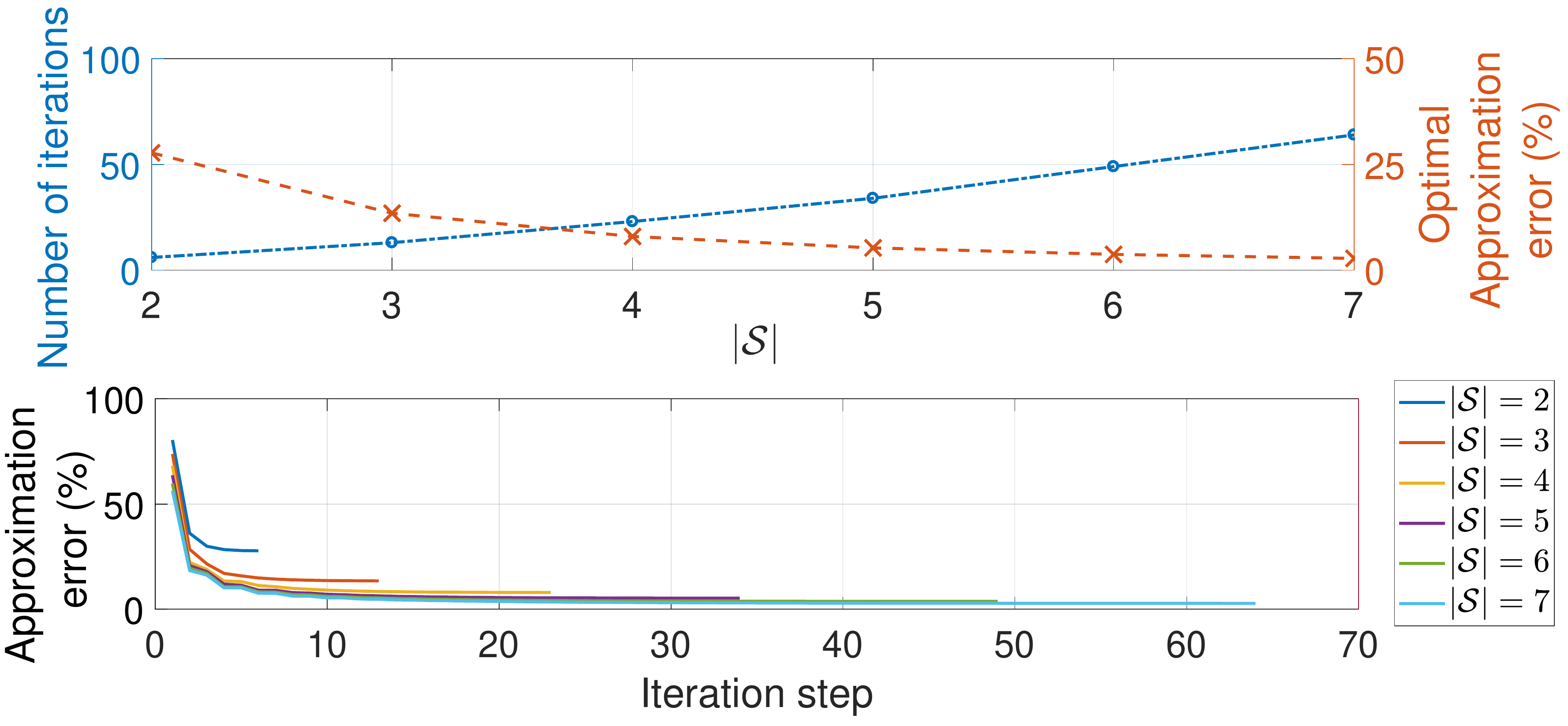}
\caption{Number of iterations and optimal approximation error as a function of $|\mathcal{S}|$ (top); Convergence of the approximation error for different $|\mathcal{S}|$ (bottom).}
\vspace{-.cm}
\label{fig: alg2}
\end{figure}

Figure~\ref{fig: alg2} shows the convergence of Algorithm 2 under different values of $|\mathcal{S}|$. We observe that as $|\mathcal{S}|$ increases, the optimal approximation error $e^{\text{max}}$ decreases and the total number of iterations grows almost linearly. 


\section{Case Studies}  \label{sec: case}

\subsection{Formulation and Setup} \label{sec: setup}
We base our chance-constrained DC OPF formulation on \cite{vrakopoulou_probabilistic_2013}; it was also used in \cite{UnimodalBL}. We assume that the system has $N_W$ wind power plants with forecast error $\tilde{w}\in \mathbb R^{N_W}$ (each element is represented as $\tilde{w}_i$), $N_G$ generators, and $N_B$ buses. The forecast errors are calculated as the difference between actual wind power realizations and their corresponding forecasts and are compensated by  reserves. Design variables include generation output $P_G \in \mathbb R^{N_G}$, up and down reserve capacities $R_G^{up} \in \mathbb R^{N_G}, R_G^{dn} \in \mathbb R^{N_G}$, and a distribution vector $d_G \in \mathbb R^{N_G}$, which determines how much reserve each generator provides to balance the total forecast error. The full problem formulation is 
\begin{subequations}
\begin{align}
\min&\ P_G^T[C_1]P_G+C_2^TP_G+C_R^T(R_G^{up}+R_G^{dn})\nonumber\\
\mathrm{s.t.}
&-P_l\leq AP_\text{inj}\leq P_l,\label{cons:line}\\
&R_G=-d_G(\Sigma_{i=1}^{N_W}\tilde{w}_i),\label{cons:resv}\\
&P_{\text{inj}}=C_G(P_G+R_G)+C_W(P_W^f+\tilde{w})-C_LP_L,\label{cons:inj}\\
&\underline{P}_G\leq P_G+R_G\leq \overline{P}_G,\label{cons:gen}\\
&-R_G^{dn}\leq R_G \leq R_G^{up},\label{cons:resvcap}\\
&\mathbf{1}_{1 \times N_G} d_G=1,\label{cons:dist}\\
&\mathbf{1}_{1 \times N_B}(C_GP_G+C_WP_W^f-C_LP_L)=0,\label{cons:bal}\\
&P_{G}\geq \mathbf{0}_{N_G \times 1},\ d_G\geq \mathbf{0}_{ N_G \times 1}, \label{cons:nneg1}\\
& R_G^{up}\geq \mathbf{0}_{N_G \times  1}, \ R_G^{dn}\geq \mathbf{0}_{N_G \times 1}, \label{cons:nneg2}
\end{align}
\end{subequations}
where $[C_1] \in \mathbb R^{N_G \times N_G}$, $C_2 \in \mathbb R^{N_G}$, and $C_R \in \mathbb R^{N_G}$ are cost parameters. Constraint \eqref{cons:line} bounds the power flows by the line limits $P_l$. The power flow is calculated from the power injections $P_\text{inj}$ in \eqref{cons:inj} and the constant matrix $A$, which is calculated using the admittance matrix and network connections. Constraint \eqref{cons:resv} computes the real-time reserve response $R_G$ that is bounded by the reserve capacities $R_G^{dn}$ and $R_G^{up}$ in \eqref{cons:resvcap}. In \eqref{cons:inj}, $P_W^f$ is the wind power forecast, $P_L$ is the load (which is assumed to be known, though the formulation can be easily extended to handle uncertain loads), and $C_G$, $C_W$, and $C_L$ are constant matrices that map generators, wind power plants, and loads to buses. Constraint \eqref{cons:gen} bounds generator outputs within their limits $[\underline{P}_G, \overline{P}_G]$; \eqref{cons:dist}, \eqref{cons:bal} enforce power balance with and without wind power forecast error; and \eqref{cons:nneg1}, \eqref{cons:nneg2} ensure all decision variables are non-negative. Uncertainty-related constraints \eqref{cons:line},\eqref{cons:gen}, and \eqref{cons:resvcap} are formulated as chance constraints. 

We test our approaches on the IEEE 118-bus and 300-bus systems modified to include a large number of wind power plants with a total of 400 and 2000 MW of forecasted wind power, respectively. We use the network and cost parameters from \cite{nesta} and set $C_R=10C_2$. We add wind power to all buses with generators and allocate the forecasted wind power to these buses in proportion their generation limit.

We also test our approaches using two forecast error data sets with different characteristics. We define the forecast error ratio as the ratio between the forecast error and the corresponding forecast. {\em Data Set 1 {\bf (DS1)}} was used in \cite{mariatsg}. The data set is generated using the Markov-Chain Monte Carlo mechanism \cite{MCMC} on real wind power forecasts and realizations from Germany. The wind power is well-forecasted with small forecast error ratios ($-30$ to $60\%$). For each wind bus, we randomly select the forecast errors from the same data pool without considering spatial correlation.
{\em Data Set 2 {\bf (DS2)}} is constructed from the RE-Europe data set \cite{reeurodata}, which contains hourly wind power forecasts and realizations based on the European energy system. The data set includes strong spatiotemporal correlation. However, the data set also contains poor forecasts with extreme forecast error ratios, up to $5300\%$ \cite{BOWENTHESIS}. Therefore, we scale down the forecast errors by $60\%$ and then filter outliers with forecast error ratios larger than $100\%$. 

We use 5000 randomly selected data points for the 118-bus system and 8000 for the 300-bus system  to construct $\mathcal{D}_{\xi}$ and $\mathcal{U}_{\xi}$. More data is needed for the 300-bus system since the uncertainty dimension is larger. In addition, we use histograms with 15 and 20 bins to determine the locations of mode $m$ for DS1 and DS2 by identifying the bin with the most points. Further, to evaluate reliability, we randomly select 5000 and 8000 data points to conduct out-of-sample tests for the 118-bus and 300-bus systems, respectively. We define the reliability as the percentage of wind power forecast errors for which all chance constraints are satisfied. To guarantee the credibility of the result, we perform three parallel tests by randomly reselecting the data used to construct the ambiguity sets.

We benchmark our approaches against two conventional approaches. 
{\em Analytical reformulation assuming multivariate Gaussian distributions} {\bf (AR)} used in \cite{bienstock_chance_2014,bowentsg,cc4}  uses  moments determined from the data. Then all chance constraints can be exactly reformulated as SOC constraints. The {\em scenario-based method} {\bf (SC)} developed in \cite{Margellos2014}  enforces the constraints affected by uncertainties to be robust against a probabilistically robust set. This set is constructed using a sufficient number of randomly selected uncertainty realizations.

We solve all optimization problems using CVX with the Mosek solver \cite{cvx1,cvx2}. We set $\epsilon=5\%$ and $\alpha=1$. The latter is valid because, in general, wind power forecast error is marginally unimodal \cite{BOWENTHESIS}.


\subsection{Comparative Results} \label{sec: result}

We first compare the DR approaches to the benchmark approaches in terms of objective cost, reliability, and computational time. The results are summarized in Table~\ref{tab: cost}. {\bf (DR-M)} is the DR approach with with ambiguity set (3), which does not include the unimodality assumption. {\bf (DR-U)} is the DR approach with ambiguity set (4), solved using the exact reformulation. To facilitate comparisons, we define a percentage difference on cost (C/Diff) and reliability (R/Diff) against the benchmarks, where AR generally produces low-cost solutions that are not sufficiently reliability and SC generally produces high-cost solutions with higher reliability than necessary. Specifically, we calculate the C/Diff of a DR approach as the difference in cost compared to that of the AR approach divided by the difference in cost between the AR and SC approaches. The R/Diff is defined similarly. Small C/Diffs are desirable, i.e., low costs approaching that of the AR approach. Large R/Diffs are desirable, i.e., high reliability approaching that of the SC approach. We define the improvement (Improv) of a DR approach to be its R/Diff divided by its C/Diff. Large Improvs are desirable, indicating a better trade-off between cost and reliability.        

\begin{table*}[t]
\setlength\tabcolsep{2pt}
\centering
\caption{Objective costs, reliability ($\%$), and computational times (seconds) for each approach}
\vspace{-.2cm}
\label{tab: cost}
\begin{tabular}{cc|rrr|rrr|rrrrrr|rrrrrr}
\hline
\multicolumn{2}{c|}{\multirow{2}{*}{Bus/Data Set}} & \multicolumn{3}{c|}{AR}     & \multicolumn{3}{c|}{SC}     & \multicolumn{6}{c|}{DR-M}                                           & \multicolumn{6}{c}{DR-U}                                              \\
\multicolumn{2}{c|}{}                              & Cost  & Reliab & \hspace{.2cm} Time  & Cost  & Reliab & \hspace{.2cm}Time  & Cost  & C/Diff & Reliab & R/Diff & Improv & \hspace{.2cm}Time  & Cost  & C/Diff & Reliab & R/Diff & Improv & \hspace{.2cm}Time    \\ \hline \hline
\multirow{3}{*}{118/DS1}           & min           & 3309  & 81.7        & 11.0  & 4935  & 100.0       & 11.2  & 3466  & 9.6       & 99.7        & 98.3             & 10.1   & 11.1  & 3340  & 1.9       & 97.0        & 83.6             & 40.4   & 475.1   \\
                                   & avg           & 3310  & 81.8        & 11.4  & 4937  & 100.0       & 14.4  & 3467  & 9.6       & 99.7        & 98.3             & 10.2   & 11.3  & 3343  & 2.0       & 97.1        & 84.2             & 41.6   & 478.4   \\
                                   & max           & 3310  & 81.9        & 11.8  & 4942  & 100.0       & 18.7  & 3468  & 9.7       & 99.7        & 98.4             & 10.2   & 11.4  & 3344  & 2.1       & 97.2        & 84.5             & 44.0   & 483.4   \\ \hline
\multirow{3}{*}{118/DS2}           & min           & 3491  & 79.5        & 11.0  & 5902  & 100.0       & 10.5  & 4064  & 23.5      & 98.9        & 94.3             & 3.3    & 11.1  & 3703  & 8.7       & 95.0        & 74.7             & 8.1    & 2490.2  \\
                                   & avg           & 3520  & 81.8        & 11.7  & 5926  & 100.0       & 11.6  & 4141  & 25.9      & 99.2        & 95.8             & 3.7    & 11.6  & 3736  & 9.0       & 95.6        & 75.7             & 8.4    & 3160.7  \\
                                   & max           & 3564  & 85.4        & 12.4  & 5942  & 100.0       & 13.3  & 4261  & 29.3      & 99.7        & 97.9             & 4.1    & 12.3  & 3780  & 9.2       & 96.6        & 76.7             & 8.7    & 3815.9  \\ \hline
\multirow{3}{*}{300/DS1}           & min           & 14408 & 72.9        & 125.4 & 18032 & 100.0       & 134.0 & 14579 & 4.7       & 99.6        & 98.5             & 20.9   & 124.5 & 14479 & 1.9       & 96.4        & 86.7             & 44.4   & 1948.5  \\
                                   & avg           & 14409 & 73.6        & 127.0 & 18038 & 100.0       & 136.9 & 14580 & 4.7       & 99.6        & 98.6             & 21.0   & 125.3 & 14479 & 1.9       & 96.6        & 87.0             & 44.9   & 2158.5  \\
                                   & max           & 14410 & 74.1        & 130.2 & 18046 & 100.0       & 140.2 & 14581 & 4.7       & 99.7        & 98.9             & 21.0   & 125.9 & 14480 & 2.0       & 96.7        & 87.4             & 45.4   & 2263.6  \\ \hline
\multirow{3}{*}{300/DS2}           & min           & 15125 & 85.7        & 234.8 & 21249 & 100.0       & 141.7 & 16956 & 29.0      & 99.7        & 97.7             & 3.1    & 236.5 & 15724 & 8.9       & 96.7        & 76.9             & 6.8    & 6547.7  \\
                                   & avg           & 15161 & 86.2        & 236.1 & 21373 & 100.0       & 143.1 & 17052 & 30.5      & 99.7        & 97.8             & 3.2    & 237.7 & 15777 & 9.9       & 97.1        & 78.8             & 8.0    & 10880.4 \\
                                   & max           & 15191 & 87.1        & 237.1 & 21436 & 100.0       & 144.9 & 17128 & 32.0      & 99.7        & 97.9             & 3.4    & 239.4 & 15880 & 11.4      & 97.4        & 81.8             & 8.7    & 15720.6 \\ \hline
\end{tabular}

\end{table*}

From Table~\ref{tab: cost}, we see that SC provides overly conservative results with the highest costs and $100\%$ reliability, AR provides the least conservative results with the lowest costs and the lowest reliability always below $95\%$, and the DR approaches provide intermediate costs and reliability, with all reliabilities above $95\%$. DR-U provides higher costs and higher reliability than DR-M since it only considers moment information in the ambiguity set. If we compare the Diffs and Improvs of DR-U and DR-M, we see that DR-U provides a better trade-off between cost and reliability. Solutions using DS1 are more stable with less variability across parallel tests than those using DS2. Additionally, solutions using DS1 have higher Improvs than those using DS2. 

\subsection{Computational Performance} \label{sec: convergence}

As shown in Table~\ref{tab: cost}, DR-U employs an iterative solution algorithm, while the other approaches do not. Hence, DR-U requires significantly more computational time. For large system dimensions, the computational burden become severe pointing to the need for approximations. Computational times are larger for DS2 than DS1.

 Table~\ref{tab: time} summarizes the percent of the total computational time to complete Step 2 of Algorithm 1 and the required number of iterations. DS2 requires a larger number of iterations than DS1 and the 118-bus system requires relatively more computational time percent to complete Step 2 than the 300-bus system. The total computation time of each iteration slightly increases
 over the iterations, while the time needed for Step 2 is approximately constant.

\begin{table}[tbp]
\caption{Algorithm 1, Step 2 percent time and number of iterations}
\vspace{-.2cm}
\centering
\label{tab: time}
\begin{tabular}{l|ccc|ccc|ccc|ccc}
\hline
\multirow{2}{*}{} & \multicolumn{3}{c|}{118/DS1} & \multicolumn{3}{c|}{118/DS2} & \multicolumn{3}{c|}{300/DS1} & \multicolumn{3}{c}{300/DS2} \\ 
                               & min      & avg     & max     & min      & avg     & max     & min      & avg     & max     & min      & avg     & max     \\ \hline\hline 
Percent                     & 86.2     & 86.8    & 87.4    & 85.4     & 85.6    & 85.9    & 59.1     & 59.1    & 59.2    & 43.1     & 43.6    & 43.9    \\
Iterations                      & 5        & 5       & 5       & 26       & 32      & 36      & 6        & 7       & 7       & 14       & 23      & 33      \\ \hline
\end{tabular}
\end{table}



\begin{figure}
\centering
\vspace{.cm}
\includegraphics[width=5in]{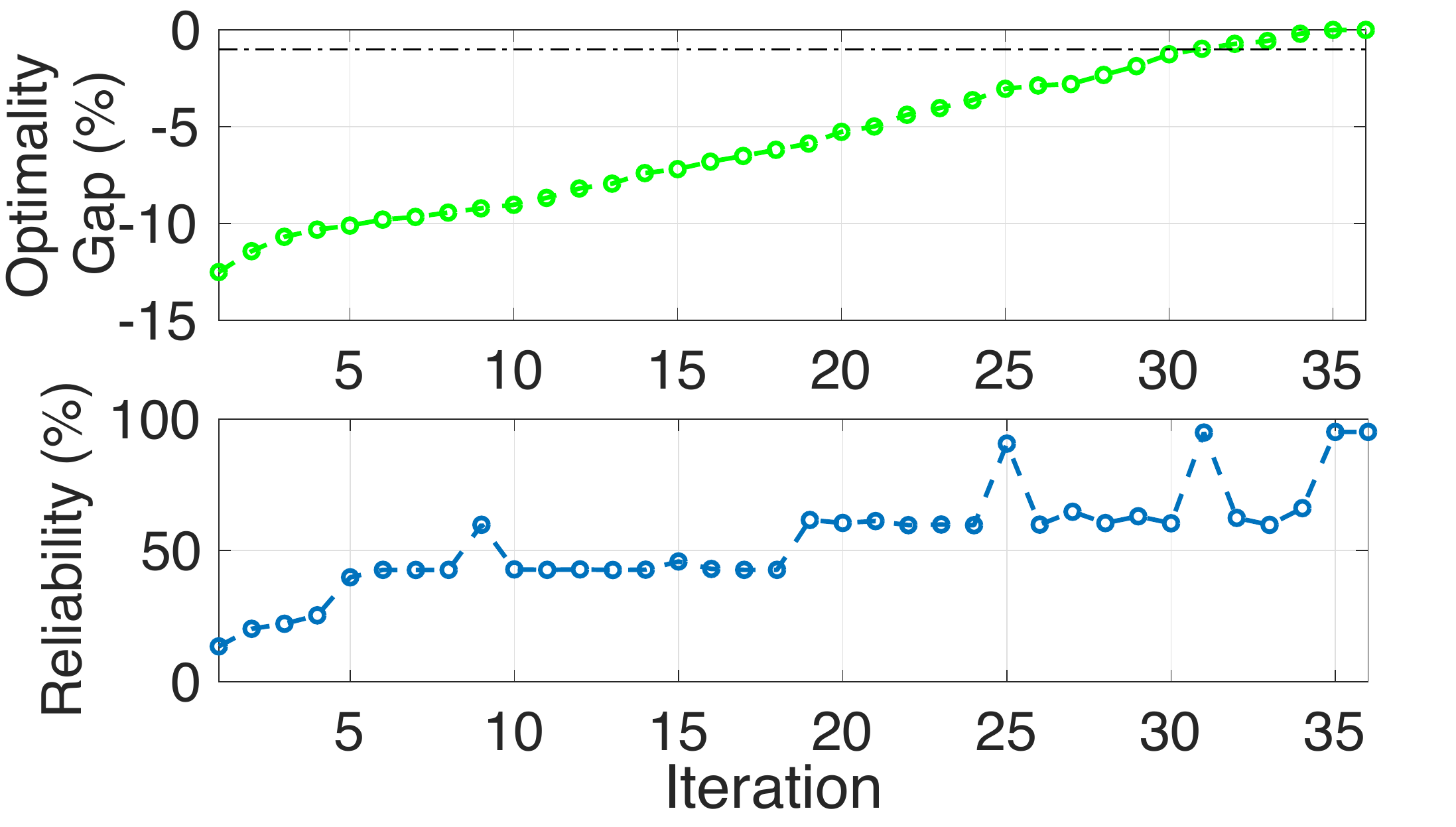}
\caption{Optimality gap and reliability of intermediate solutions of Algorithm~1. Black dashed line marks $1\%$ optimality gap.}
\vspace{-.cm}
\label{fig: lbinner}
\end{figure}

Next, we check if the solutions from the intermediate iterations of Algorithm 1 are good approximates of the  optimal solution. Fig.~\ref{fig: lbinner} shows the optimality gap and reliability of the intermediate solutions for the 118-bus system using DS2. Note that each intermediate solution constitutes a relaxed approximation and so the optimality gap is negative. We find that the intermediate solutions are not good approximates because solutions with small absolute optimality gaps ($<1\%$) can have low reliability ($<70\%$). We also observe that higher objective cost does not always guarantee higher reliability in out-of-sample tests.

\subsection{Conservative Approximations}
In this section, we compare five options for generating conservative approximations based on Proposition~\ref{prop: drcc_set2_ub}:
\begin{itemize}
\item UB: the approach in \cite{UnimodalBL} without OPS. 
\item OPS0: an online approach that uses the OPS solutions with $|\mathcal{S}|=K-1$ on the violated constraints from Step~2 in Algorithm 1.
\item OPS1: an offline approach that uses the OPS solutions with $|\mathcal{S}|=K-1$ on all of the DR chance constraints.
\item OPS2: an aggregated version of OPS0 that uses all the OPS solutions with $1\leq|\mathcal{S}|\leq K-1$.
\item OPS3: an aggregated version of OPS1 that uses all the OPS solutions with $1\leq|\mathcal{S}|\leq K-1$.
\end{itemize}
The online options (OPS0 and OPS2) require information about which DR chance constraints are violated in Step~2 of Algorithm~1, but the offline options (OPS1 and OPS3) simply apply the approximations to all DR chance constraints. Further, the aggregated versions (OPS2 and OPS3) take advantage of OPS solutions with smaller parameter dimension. 

\begin{figure*} 
    \centering
  \subfloat[118/DS1 \label{fig: 118ds1}]{%
       \includegraphics[width=0.48\linewidth]{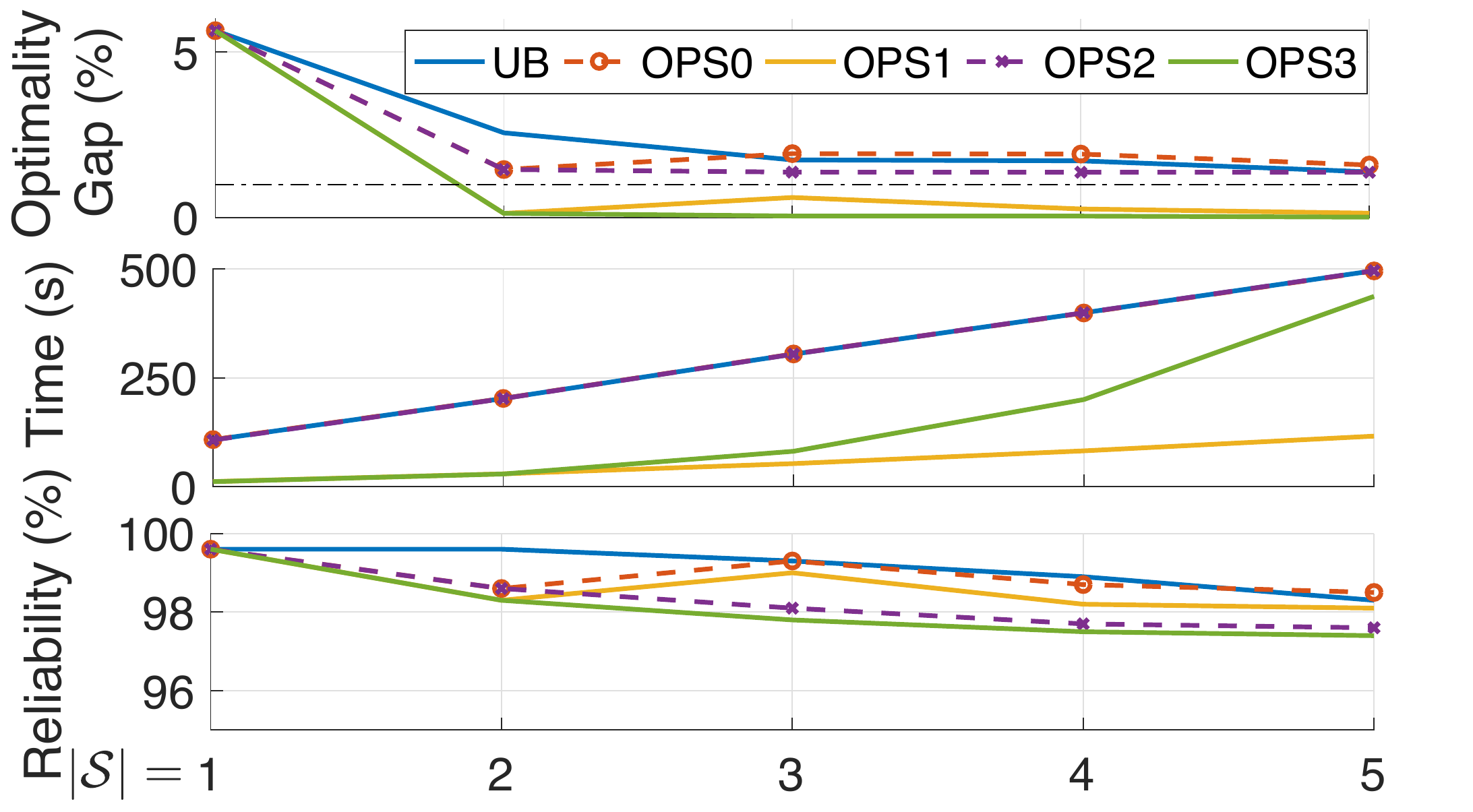}}
    \hfill
  \subfloat[300/DS1 \label{fig: 300ds1}]{%
        \includegraphics[width=0.48\linewidth]{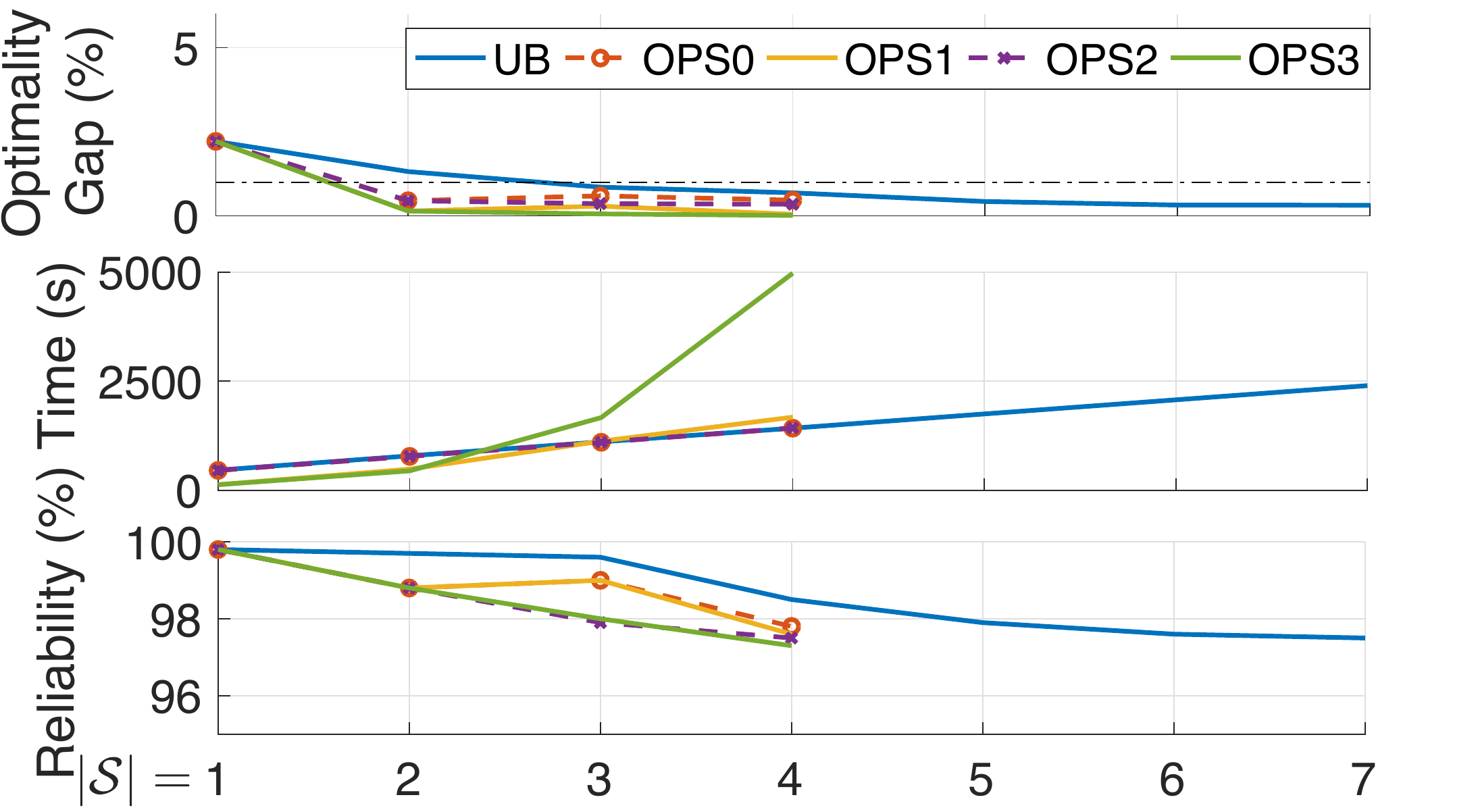}}
    \\
  \subfloat[118/DS2 \label{fig: 118ds2}]{%
        \includegraphics[width=0.48\linewidth]{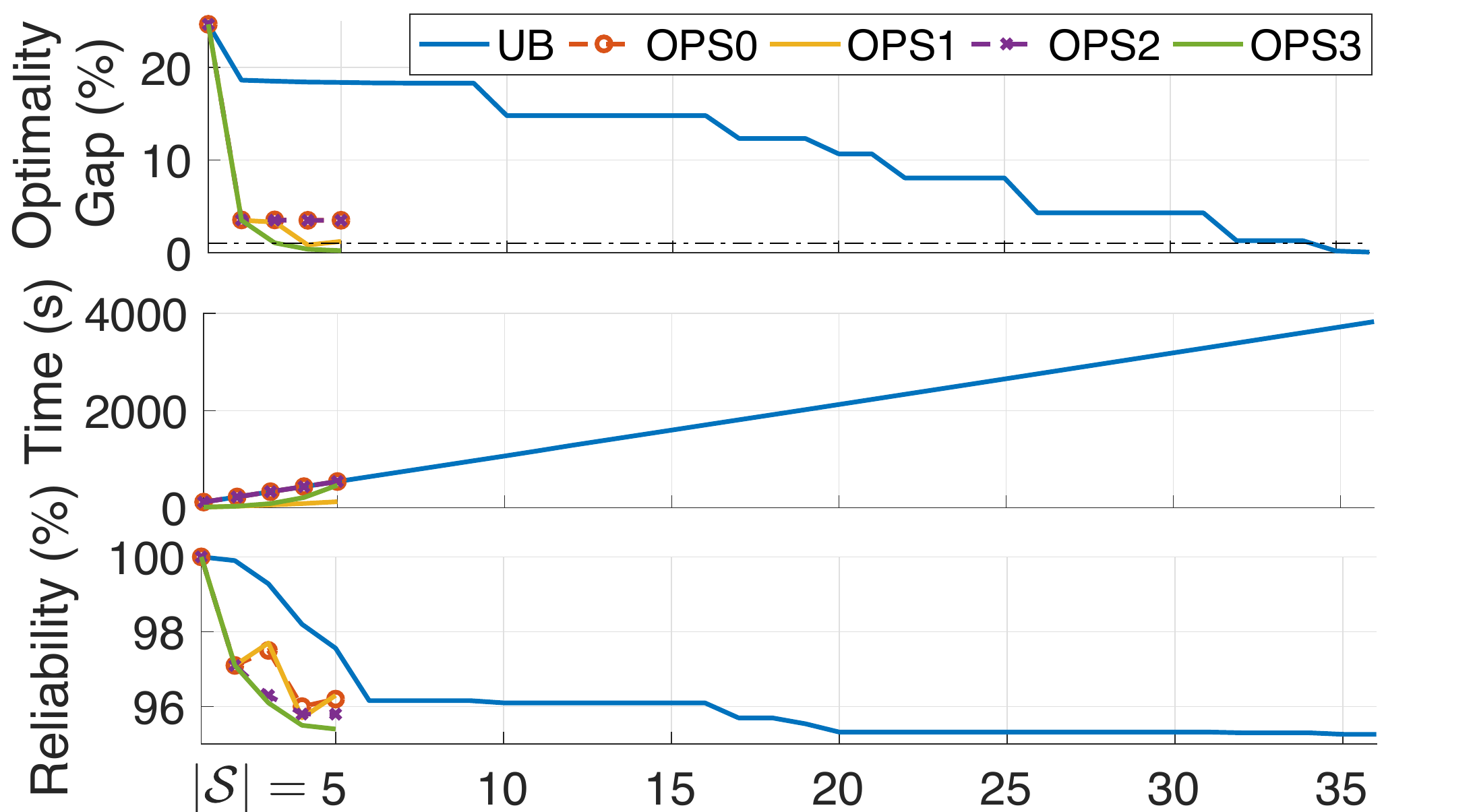}}
    \hfill
  \subfloat[300/DS2 \label{fig: 300ds2}]{%
        \includegraphics[width=0.48\linewidth]{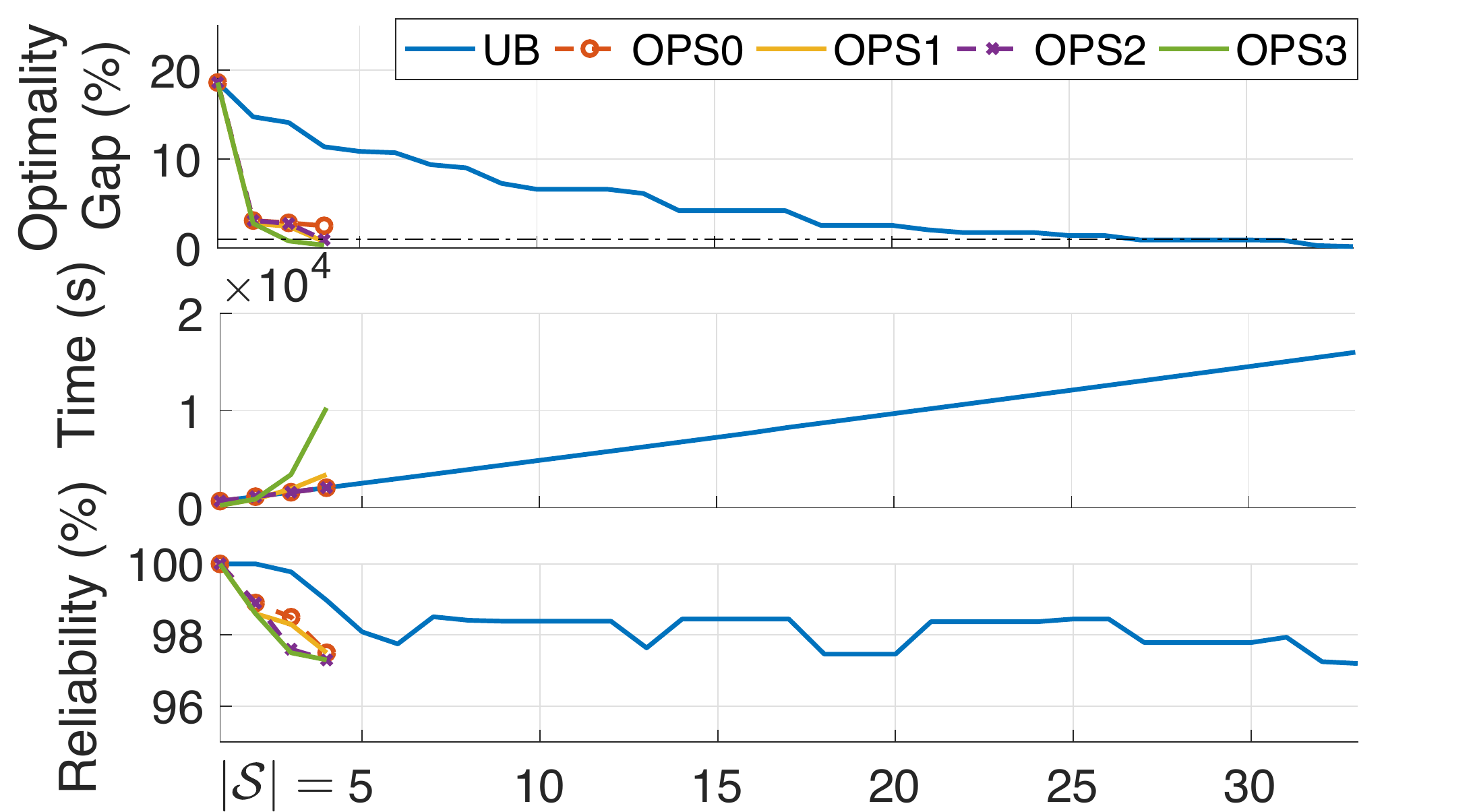}}
  \caption{Performance of the conservative approximations. Black dashed lines in top subplots mark $1\%$ optimality gap. Note that, in some cases, comparable plots use different scales so that the difference between the OPS options is visible.}
  \label{fig:approx} 
\end{figure*}

Figure~\ref{fig:approx} compares all of the approximations on both test systems using both data sets. For the OPS options we limit $|\mathcal{S}|\leq 5$ for the 118-bus system and $|\mathcal{S}|\leq 4$ for the 300-bus system. For the 118 bus system with DS1, Fig.~\ref{fig: 118ds1} shows that UB fails to achieve a $\leq 1\%$ optimality gap. OPS0 and OPS1 demonstrate better convergence rates and optimality gaps but their optimality gaps (i.e., costs) do not continue to decrease as $|\mathcal{S}|$ increases. On the other hand, by taking advantage of OPS solutions with smaller $|\mathcal{S}|$, 
OPS2 and OPS3 have similar convergence rates as well as non-increasing optimality gaps as $|\mathcal{S}|$ increases. All approximation options (except UB, OPS0, OPS1 for $K-1=5$) take less time than  exactly solving DR-U (483.4s). The offline options (OPS1 and OPS3) take a similar amount of time as AR, SC, and DR-M when $K-1$ is small, while the online options exhibit a linear relationships between computational time and $K-1$. All approximate solutions satisfy the $95\%$ constraint satisfaction level. Tighter approximations (i.e., larger $K-1$) are less conservative leading to lower reliability.

For the 300-bus system with DS1, Fig.~\ref{fig: 300ds1} shows similar trends except that the computational times of OPS1 and OPS3 become larger than those of the online options when $|\mathcal{S}|\geq 3$. When $|\mathcal{S}|=2$, solutions from all OPS options achieve $\leq 1\%$ optimality gaps, $\geq 95\%$ reliability, and computational times less than exactly solving DR-U (2263.6s). 

Figures~\ref{fig: 118ds2} and~\ref{fig: 300ds2} show that the approximations with DS2 generally have larger optimality gaps than those with DS1. Also, with DS2, UB requires many parameters while all OPS options converge with much fewer parameters and with computational times less than exactly solving DR-U (3815.9s for the 118 bus system and 15720.6s for the 300 bus system). In Fig.~\ref{fig: 300ds2} we again observe that offline OPS options can take more computational time than the online options as $|\mathcal{S}|$ increases. In general, the offline options are less computationally advantageous for the 300 bus system than the 118 bus system. Further, in the reliability plot in Fig.~\ref{fig: 300ds2} we observe oscillations like in Fig.~\ref{fig: lbinner} demonstrating, again, that reliability does not always increase with higher objective costs.

In summary, the conservative approximations produce good approximates of the optimal solution of DR-U with small optimality gaps and high reliability. Further, the OPS options improve upon the previously-developed approach (UB) by achieving much better convergence rates and solution quality.

\section{Conclusions} \label{sec: conclusion}
In this paper, we developed an OPS approach to achieve a better conservative approximation of a DR chance constraint with moment and unimodality information. We further proposed multiple online and offline options to generate the approximation and evaluated their performance against the current state of the art. Through case studies on modified IEEE 118-bus and 300-bus systems, we demonstrated that including unimodality information within a DR OPF problem with wind power uncertainty leads to a better cost/reliability trade-off than benchmark approaches or a DR OPF that includes only moment information. However, it also increased the computational time.  We showed that conservative approximations reduce the computational time and options leveraging our OPS approach provide solutions with low optimality gaps, that satisfy desired reliability levels, and that require less computational time. We also showed how the results vary across two forecast error data sets. We found that both the data set and choice of test system have a significant impact on the value of including unimodality information in DR OPF, indicating that, in practice, the value is highly system-dependent. Moreover, the relative performance of algorithmic approaches, in terms of optimality gap, computational time, and solution reliability, is also system-dependent.




\appendices
\section{Proof of Theorem~\ref{thm: pwl_multi}}\label{appd: proof1}
Condition 1: The last piece of $h(\tau)$, i.e., $h_{|S|}(\tau)$, must have zero slope because otherwise the error is infinite (if the slope is strictly positive) or $h_{|S|}(\tau) < v(\tau)$ for a sufficiently large $\tau$ (if the slope is strictly negative). It follows that $h_{|S|}(\tau) = \lim_{\tau \rightarrow \infty} v(\tau) = \sqrt{(1-\epsilon)/\epsilon}$ because this is the constant function that dominates $v(\tau)$ with the smallest error.

Condition 2: Since $v(\tau)$ is non-decreasing and concave, we have $d_s \geq 0$ and $d_s$ is non-increasing in $s$. If $h_s(\tau)$ is not tangent to $v(\tau)$, then we decrease $f_s$ until $h_s(\tau)$ is tangent to $v(\tau)$. Note that this does not increase $e^{\text{max}}$. Then, all pieces of $h(\tau)$ are tangent to $v(\tau)$, except $h_1(\tau_0) = v(\tau_0)$ where $h'_1(\tau_0) > v'(\tau_0)$. In this case, we rotate $h_1(\tau)$ clockwise around the point $(\tau_0, v(\tau_0))$ until $h_1(\tau)$ becomes tangent to $v(\tau)$ at $\tau_0$. Note that the rotation does not increase $e^{\text{max}}$.

Condition 3: We prove by contradiction. Assume that $h^t(\tau)$ satisfies all three conditions and has an error $e^{t,\text{max}}$, and there exists an $h^c(\tau)$ that satisfies Conditions 1 and 2 and has an error $e^{c,\text{max}} < e^{t,\text{max}}$. If $|\mathcal{S}|=1$, then $h^t(\tau) = h^c(\tau)$ due to Condition 1. This contradicts the assumption. If $|\mathcal{S}|>1$, since $e^{c,\text{max}}<e^{t,\text{max}}=E^t_2=E^t_{|\mathcal{S}|}$, we have $B_2^t>B_2^c$ and $B_{|\mathcal{S}|}^t<B_{|\mathcal{S}|}^c$. If $|\mathcal{S}|=2$, this is a clear contradiction. If $|\mathcal{S}|>2$, then there exists an $s \in [2, |S|-1]$ such that $[B_s^t,B_{s+1}^t]\subsetneqq [B_s^c,B_{s+1}^c]$, i.e., there exists a pair of pieces $h^t_s(\tau)$ and $h^c_s(\tau)$ with the same index $s$ such that the domain of $h^t_s(\tau)$ is a strict subset of that of $h^c_s(\tau)$, because $h^t(\tau)$ and $h^c(\tau)$ have the same domain $[\tau_0, \infty)$ and the same number of pieces. According to Condition 2, both $h^t_s(\tau)$ and $h^c_s(\tau)$ are tangent to $v(\tau)$ and hence $e^{t,\text{max}} = E^t_s = E^t_{s+1} \leq \max\{E^c_s,E^c_{s+1}\} \leq e^{c,\text{max}}$, contradicting the assumption.

\section{Proof of Theorem~\ref{thm: pwl_multi_exist}}\label{appd: proof2}
First, with a similar proof to that of Theorem~\ref{thm: pwl_multi}, we can show that an $|S|$-piece PWL function $h(\tau)$ is an optimal outer approximation of $v(\tau)$ on a bounded interval $[\tau_0, \overline{\tau}]$ if it satisfies Conditions 2 and 3 in Theorem~\ref{thm: pwl_multi}. We term this result Theorem~\ref{thm: pwl_multi}-finite. 

Second, we use mathematical induction to prove that there exists an $|\mathcal{S}|$-piece PWL approximation that satisfies the conditions of Theorem~\ref{thm: pwl_multi}-finite, if $v(\tau)$ is defined on a bounded interval $[\tau_0, \overline{\tau}]$. When $|\mathcal{S}|=1$, Condition 3 becomes $h(\tau_0)-v(\tau_0)=h(\overline{\tau})-v(\overline{\tau})$. The single-piece optimal PWL approximation exists by simply searching for a point in $[\tau_0,\overline{\tau}]$ at which $h(\tau)$ and $v(\tau)$ are tangent.

Finally, we show that if a $C$-piece optimal PWL approximation exists and satisfies the conditions of Theorem~\ref{thm: pwl_multi}-finite, then so does a $(C+1)$-piece optimal PWL approximation. In the $C$-piece approximation, denote the second largest break point as $B_F$.  Consider a $C$-piece approximation on $[\tau_0,B_F]$ and a single-piece approximation on $[B_F,\overline{\tau}]$. As we move $B_F$ from $\tau_0$ to $\overline{\tau}$, the error of the $C$-piece approximation continuously increases from zero to a finite positive number (i.e., the optimal error for a $C$-piece approximation on $[\tau_0,\overline{\tau}]$) while the error of the single-piece approximation continuously decreases from a finite positive value (i.e., the optimal error for a single-piece approximation on $[\tau_0,\overline{\tau}]$) to zero. It follows that there exists a $B^*_F \in [\tau_0, \overline{\tau}]$ such that the error of the $C$-piece approximation (on $[\tau_0, B^*_F$]) equals that of the single-piece approximation (on $[B^*_F, \overline{\tau}]$). The resultant $(C+1)$-piece PWL approximation satisfies the conditions of Theorem~\ref{thm: pwl_multi}-finite. The proof of the case with $\overline{\tau} = +\infty$ is similar and so omitted.

\bibliography{ReferencesJournal1}
\bibliographystyle{IEEEtran}

\end{document}